# New Residue Arithmetic Based Barrett Algorithms, Part II: Modular Polynomial Computations

Hari Krishna Garg[1] & Hanshen Xiao[2]

*Abstract*: **In this paper, we derive a new computational algorithm for Barrett technique for modular polynomial multiplication, termed BA-P. BA-P is then applied to a new residue arithmetic based Barrett algorithm for modular polynomial multiplication (BA-MPM). The focus of the work is an algorithm that carries out the entire computation using only modular arithmetic without conversion to large degree polynomials. There are several parts to this work. First, we set up a new BA-P using polynomials other than $u^{\alpha}$. Second, residue arithmetic based BA-MPM is described. A complete mathematical framework is described including proofs of the steps in the computations and the validity of results. Third, we present a computational procedure for BA-MPM. Fourth, the BA-MPM is used as a basis for algorithms for modular polynomial exponentiation (MPE). Applications are in areas of signal security and cryptography.**

**Keywords:** Barrett Algorithm (BA), BA for Polynomials (BA-P), Modular Polynomial Multiplication (BA-MPM), Montgomery Multiplication (MM), Residue Polynomial Systems (RPS), Chinese Remainder Theorem for Polynomials (CRT-P), Modular Polynomial Multiplication (MPM), Modular Polynomial Exponentiation (MPE), Base Extension for Polynomials (BEX-P), Modular Addition (MADD), Modular multiplications (MMULT), Mixed Radix System (MRS).

[1]Hari Krishna Garg (eleghk@nus.edu.sg), is with the Electrical & Computer Engineering Department, National University of Singapore, Singapore.
[2]Hanshen Xiao (xhs13@mails.tsinghua.edu.cn) is with the Department of Mathematics, Tsinghua University, Beijing, China.

The work of Hari K Garg was carried out at the NUS-ZJU Sensor-Enhanced Social Media (SeSaMe) Centre. It is supported by the Singapore National Research Foundation under its International Research Centre @ Singapore Funding Initiative and administered by the Interactive Digital Media Program Office.

The work of Hanshen Xiao was supported by Tsinghua University Initiative Scientific Research Program 20141081231, the National Natural Science Foundation of China under Grant 61371104, and National High Technology Project of China under Grant 2011AA010201.

## I. Introduction

CRYPTOGRAPHY TECHNIQUES play an important role in the security of communication and data processing systems. Instances of such cryptography techniques include RSA (Rivest-Shamir-Adelman), Rabin, Diffie Hellman and El Gamal. Many of these and other related techniques deal with arithmetic defined in a large size finite fields GF($p$) and/or GF($p^N$), where $p$ is a prime integer and $N$ is any positive integer. A large size field may be realized by setting $p = 2$ (binary arithmetic) and $N$ to be a large value in order to provide adequate level of security. This paper also focuses on such finite fields with values of $N$ running from several hundreds to several thousands. The elements in GF($2^N$), are represented as polynomials of degree up to $N - 1$. A challenge in such environments is to perform the arithmetic in a computationally efficient and timely manner. The basic arithmetic operations are:

1. Multiplication of two elements in GF($2^N$):

   $C(u) = A(u) \cdot B(u) \pmod{P(u)}; 0 \le \deg(A(u)), \deg(B(u)) < \deg(P(u)) = N$; and

2. Modular exponentiation in GF($2^N$):

   $C(u) = A(u)^E \pmod{P(u)}; 0 \le \deg(A(u)) < \deg(P(u)) = N$.

Here, $P(u)$ is a monic irreducible (or primitive) polynomial in GF(2). This eliminates a straightforward algorithm to compute $C(u)$ based on the Chinese remainder theorem for polynomials (CRT-P). Ordinary computations, 1 and 2 as above without the mod $P(u)$ operation, are relatively simpler to perform. However mod $P(u)$ operation requires polynomial division and it is a difficult operation to carry out. Further, modular polynomial exponentiation (MPE) is carried out via repeated use of modular polynomial multiplication (MPM) algorithm. Hence it is essential that an efficient algorithm be used in carrying out MPM. It is observed here that the arithmetic that may appear to be rather straightforward may turn out to be a computational bottleneck in applications such as cryptography due to involvement of large degree polynomials. In such applications, degree of polynomials $A(u)$, $B(u)$, and $P(u)$ may be several hundred to a few thousand in value.

We reiterate that the task that we are dealing with is that of polynomial multiplication carried out mod $P(u)$.

Hence the results are general and may find applications beyond finite field arithmetic and cryptography.

Residue arithmetic is extensively used to express a large size computation in a polynomial ring into a number of smaller size computations in several polynomial rings. Computations in these smaller rings can be carried out in parallel. Residue number systems have been applied extensively in combination with computational methods, such as Barrett algorithm (BA) and Montgomery multiplication (MM), to carry out modulo arithmetic operations in large size finite integer rings. Residue Polynomial Systems (RPS) have also been used extensively for computing polynomial multiplication in digital signal processing. However there is a distinct gap when it comes to RPS based BA and MM for polynomials.

The major contributions of the paper are as follows. The primary and eventual objective of the work is to compute MPM and MPE as efficiently as possible with a view towards their application to signal security and cryptography. This is accomplished by first describing a new BA for modular polynomial multiplication (BA-P) for computing the quotient polynomial $C(u)$ associated with a polynomial $X(u)$ when it is divided by a polynomial $P(u)$. It is assumed that $N = \deg(P(u))$ is a large integer. Second, a residue arithmetic based BA-P, termed BA-MPM, is described for modular polynomial multiplication. Third, a computationally efficient procedure for the computations underlying the new BA-MPM is described. Fourth, the new BA-MPM is used as a basis for a new algorithms for MPE. The results as general though described for GF(2) in many instances. They are valid for polynomials defined in other fields such as GF($p^N$), $N \geq 1$, and rational, real, and complex numbers.

Mathematically, RPS represents a polynomial ring modulo $M(u) = \prod_{i=1}^{n} M_i(u)$, $L = \deg(M_1(u))$, as a direct sum of $n$ smaller degree polynomial rings modulo $M_i(u)$, $L_i = \deg(M_i(u))$, $i = 1, 2, \ldots, n$. Modulo will be written as 'mod' from here onwards. The factors of $M(u)$ are pair-wise co-prime, that is, $\gcd(M_i(u), M_j(u)) = 1$, $i \neq j$. Also, $L = \sum_{i=1}^{n} L_i$ . These $n$ polynomial rings can operate independently and hence in parallel to each other. It is desirable, therefore, to stay in the smaller degree rings as much as possible and not carry out the computation back to the larger degree ring. This is even more urgent in cryptography when the polynomials can have large degrees. This is also our secondary objective.

There is abundance of research materials on Montgomery and Barrett multiplication algorithms and their various extensions. A partial list of such papers is [1]-[13]. The polynomials versions of Montgomery and Barrett multiplication algorithms can also be found in many places [14]-[22]. In [15] the authors present a digit-serial multiplication in GF($2^N$) based on Barrett modular reduction. An improved version of digit-serial multiplication algorithm is described in [16]. Other aspects such as avoiding the pre-computation phase have also been explored [19].

However, there is no paper thus far on using the richness of CRT-P and RPS to carry out MPM and MPE computation via BA. It is this particular aspect that we deal with in this paper. The algorithm described here begins with a reformulation of the basic approach to BA. It then applies RPS in a way that computations for the new BA-P stay within the residue arithmetic.

The organization of this paper is as follows. Section II provides mathematical preliminaries on polynomial arithmetic, BA-P, RPS, CRT-P, mixed radix systems for polynomials (MRS-P), and base extension for polynomials (BEX-P). Sections III-V are on the key contributions of the work. Mathematical structure for the new BA-P is described in Section III. The computational steps for the corresponding RPS based BA-P algorithm are presented in Section IV. Examples are presented to illustrate the algorithm. In Section V, we describe an algorithm for MPE that uses the new BA-MPM. Section VI is on conclusions of the work.

## II. Mathematical Preliminaries

**Polynomial Arithmetic.** Given two polynomials $X(u)$ and $A_1(u)$ with coefficients defined in a field **F**, consider dividing $X(u)$ by $A_1(u)$ and writing

$$X(u) = Q_1(u) \cdot A_1(u) + R_1(u), \qquad (1)$$

where $Q_1(u)$ is the quotient, $R_1(u)$ is the remainder, and $0 \leq \deg(R_1(u)) < \deg(A_1(u))$. Also, $\deg(Q_1(u)) = \deg(X(u)) - \deg(A_1(u))$ with $Q_1(u) = 0$ if $\deg(X(u)) < \deg(A_1(u))$. We also write (1) as

$$X(u) \equiv R_1(u) \pmod{A_1(u)}. \quad (2)$$

Dividing both sides of (1) by $A_1(u)$, we get,

$$X(u) / A_1(u) = Q_1(u) + R_1(u) / A_1(u). \quad (3)$$

The last term on the right when expressed as a sum of powers of $u$ will only contain negative powers. We also write

$$Q_1(u) = \lfloor X(u) / A_1(u) \rfloor, \quad (4)$$

$\lfloor Y(u) \rfloor$ being the polynomial floor function of $Y(u)$. We observe that $Q_1(u)$ and $R_1(u)$ are unique. The process can be repeated between $Q_1(u)$ and another polynomial, say $A_2(u)$. Thus

$$Q_1(u) = Q_2(u) \cdot A_2(u) + R_2(u), \quad (5)$$

$0 \leq \deg(R_2(u)) < \deg(A_2(u))$. To take this one step further, replace $Q_1(u)$ in (3) by the expression in (5) to get,

$$X(u) = Q_2(u) \cdot [A_2(u) \cdot A_1(u)] + [R_2(u) \cdot A_1(u) + R_1(u)]. \quad (6)$$

A generalization of (6) leads to

$$X(u) = Q_a(u) \cdot [(A_a(u) \cdots A_1(u)] + [R_a(u) \cdot \{A_{a-1}(u) \cdots A_1(u)\} + \cdots + R_2(u) \cdot A_1(u) + R_1(u)]. \quad (7)$$

This expression would be most useful in various computations involving RPS and BEX-P.

**Barrett Algorithm for Polynomials (BA-P).** Given polynomials $A(u)$, $B(u)$, and $P(u)$, $0 \leq \deg(A(u)$, $\deg(B(u) < \deg(P(u)) = N$, and the task of computing MPM

$$C(u) = A(u) \cdot B(u) \pmod{P(u)}, \quad (8)$$

BA-P, consists in computing the quotient $Q(u)$ first such that

$$X(u) = A(u) \cdot B(u) \quad (9)$$

$$X(u) = Q(u) \cdot P(u) + C(u);\ 0 \le \deg(C(u)) < N. \tag{10}$$

Then $C(u)$ is computed as

$$C(u) = X(u) - Q(u) \cdot P(u). \tag{11}$$

The computations in (9) and (11) are additional to BA required to carry out MPM. The BA-P is the computation of $Q(u)$ from $X(u)$ and $P(u)$ by expressing $Q(u)$ in (11) as

$$Q(u) = \lfloor X(u) / P(u) \rfloor, \tag{12}$$

In the current versions of BA-P [14]-[22], the computation in (12) is expressed as

$$Q(u) = \lfloor X(u) / u^a \cdot \Omega(u) / u^{a+b} \rfloor \ \approx \lfloor \lfloor X(u) / u^a \rfloor \cdot \lfloor \Omega(u) \rfloor / u^b \rfloor, \tag{13}$$

where $\lfloor \Omega(u) \rfloor$ is a pre-computed polynomial given by

$$\mu(u) = \lfloor \Omega(u) \rfloor = \lfloor u^{a+b} / P(u) \rfloor.$$

The values of scalars $a$ and $b$ can be chosen such that $Q(u)$ as computed in (13) is same as $Q(u)$ in (10) [15, 16]. The BA-P consists in the following steps:

0. Pre-Compute $\mu(u) = \lfloor u^{a+b} / P(u) \rfloor$
1. Compute $X(u) = A(u) \cdot B(u)$
2. Compute $D(u) = \lfloor X(u) / u^a \rfloor$
3. Compute $E(u) = D(u) \cdot \mu(u)$
4. Compute $Q(u) = \lfloor E(u) / u^b \rfloor$
5. Compute $C(u) = X(u) - Q(u) \cdot P(u)$.

The readers are referred to [23]-[28] for a fuller description of BA and MM, and their applications to cryptography.

**Residue Polynomial System (RPS) [23, 29].** A RPS is a polynomial ring defined by $n$ pairwise co-prime *monic* polynomials $M_1(u)$, $M_2(u)$, ... , $M_n(u)$ with elements in a field **F**. The elements of such a polynomial ring are polynomials of degree up to $L - 1$, $L = \deg(M(u))$, where

$$M(u) = \prod_{i=1}^{n} M_i(u). \tag{14}$$

Arithmetic operations [+, ·] in the polynomial ring are defined mod $M(u)$. A polynomial $X(u)$ in the given RPS is represented in terms of $n$ residue polynomials,

$$X(u) \leftrightarrow \mathbf{X}(\boldsymbol{u}) \leftrightarrow [X_1(u)\, X_2(u) \cdots X_n(u)], \tag{15}$$

where

$$X_i(u) \equiv X(u) \pmod{M_i(u)},\ i = 1, 2, ... , n. \tag{16}$$

RPS are used to express a large degree polynomial computation into a direct sum of smaller degree polynomial computations that can be carried out in parallel. For instance, the polynomial product $C(u) = A(u) \cdot B(u) \bmod M(u)$ can be computed as $n$ parallel products $C_i(u) = A_i(u) \cdot B_i(u) \pmod{M_i(u)}$, $[A_i(u), B_i(u), C_i(u)] \equiv [A(u), B(u), C(u)] \pmod{M_i(u)}$, $i = 1, \ldots , n$.

**Chinese remainder theorem for polynomials (CRT-P) [2, 3, 23, 29].** Given $\mathbf{X}(\boldsymbol{u})$, computation of $X(u)$, $0 \leq \deg(X(u)) < L$, can be done via CRT-P, stated as

$$X(u) \equiv \sum_{i=1}^{n} T_i(u) \cdot X_i(u) (\bmod M_i(u)) \cdot \left( \frac{M(u)}{M_i(u)} \right). \tag{17}$$

The reconstruction polynomials $T_i(u)$, $\deg(T_i(u)) < \deg(M_i(u))$, are computed a-priori by solving the congruence

$$T_i(u) \cdot \left( \frac{M(u)}{M_i(u)} \right) \equiv 1 \left(\bmod M_i(u)\right),\ i = 1, 2, \ldots , n. \tag{18}$$

It is clear from (18) that

$$\gcd(T_i(u), M_i(u)) = 1, i = 1, 2, \dots, n. \tag{19}$$

The CRT computation of $X(u)$ from $\mathbf{X}(\boldsymbol{u})$ may involve large degree polynomials as the dynamic range of RPS may be large; especially in applications pertaining to cryptography.

**Mixed Radix System (MRS).** MRS can also be used to compute $X(u)$ from its residues $\mathbf{X}(\boldsymbol{u})$. $X$ is represented in the following manner in an MRS:

$$X(u) = Y_1(u) + Y_2(u) \cdot M_1(u) + \cdots + Y_n(u) \cdot (M_1(u) \cdot M_2(u) \cdots M_{n-1}(u)). \tag{20}$$

The mixed radix polynomials satisfy $0 \le \deg(Y_i(u)) < \deg(M_i(u))$. $Y_i(u)$, $i = 1, \dots, n$, can be computed from $\mathbf{X}(\boldsymbol{u})$ using RPS to MRS conversion algorithms.

**Base Extension for Polynomials (BEX-P).** Consider $\mathbf{X}(\boldsymbol{u})$, residues of a polynomial $X(u)$ in (15). In many situations, one requires computation of $t$ additional residues of $X(u)$, namely $X_j(u) \equiv X(u) \pmod{M_j(u)}$, $j = n + 1, \dots, n + t$, in yet another RPS defined by $t$ additional moduli

$$M_{\mathrm{I}}(u) = \prod_{j=n+1}^{n+t} M_j(u), \tag{21}$$

where $\gcd(M(u), M_{\mathrm{I}}(u)) = 1$. This is called BEX-P. BEX-P tends to be expensive computationally. Using (17), we can compute BEX-P as follows:

$$X_j(u) \equiv \left\{ \sum_{i=1}^{n} \{T_i(u) \cdot X_i(u) (\operatorname{mod} M_i(u))\} (\operatorname{mod} M_j(u)) \cdot \frac{M(u)}{M_i(u)} (\operatorname{mod} M_j(u)) \right\} (\operatorname{mod} M_j(u)), \tag{22}$$

$j = n + 1, \dots, n + t$ [22]. Given $X_i(u)$ and pre-computed quantities such as $T_i(u)$ and $\frac{M(u)}{M_i(u)} (\operatorname{mod} M_j(u))$, $i = 1, \dots, n$, computation of $x_j(u)$, $j = n + 1, \dots, n + t$, in (22) requires $(n - 1) \cdot t$ modular additions (MADD), $n \cdot (t + 1)$ modular multiplications (MMULT), and $n \cdot t$ modular reductions. Another algorithm for BEX-P is presented in Appendix B.

## III. A New Barrett Algorithm for Polynomials

A RPS based Montgomery multiplication algorithm has been described in [22]. However, there is no such algorithm for the BA-P. We cite [15, 16, 19] and the references therein. The mod polynomials that have been used are of the type $u^a$ and hence they don't lend themselves to RPS.

Here, we first revisit the computation of $Q(u)$ in (12). We now introduce two polynomials $G(u)$ and $H(u)$, not necessarily of the form $u^a$, and approximate $Q(u)$ in (12) as

$$Q(u) = \lfloor X(u) / G(u) \cdot \Omega(u) / H(u) \rfloor \approx \lfloor \lfloor X(u) / G(u) \rfloor \cdot \mu(u) / H(u) \rfloor, \tag{23}$$

where $\mu(u)$ is a pre-computed polynomial given by

$$\mu(u) = \lfloor \Omega(u) \rfloor = \lfloor G(u) \cdot H(u) / P(u) \rfloor. \tag{24}$$

Since $X(u)$ is obtained by taking the ordinary product of $A(u)$ and $B(u)$, $0 \leq \deg(A(u)), \deg(B(u)) < N$, we have $\deg(X(u)) \leq 2 \cdot N - 2$.

Now we derive conditions on polynomials $G(u)$ and $H(u)$ for the approximation of $Q(u)$ in (23) to be equal to actual $Q(u)$ in (12). $\lfloor V(u) \rfloor$ is the polynomial part of $V(u)$. Consider dividing $V(u)$ by $S(u)$ to write $V(u) = Q(u) \cdot S(u) + R(u)$. Then we have

$$V(u) / S(u) = Q(u) + R(u) / S(u) = \lfloor V(u) / S(u) \rfloor + \delta(u), \deg(\delta(u)) \leq -1. \tag{25}$$

Applying (25) to (23), we get

$$Q(u) = \left\lfloor \frac{\left( \left\lfloor \frac{X(u)}{G(u)} \right\rfloor + \delta(u) \right) \cdot \left( \left\lfloor \frac{G(u) \cdot H(u)}{P(u)} \right\rfloor + \phi(u) \right)}{H(u)} \right\rfloor$$

$$= \left\lfloor \frac{\left\lfloor \frac{X(u)}{G(u)} \right\rfloor \cdot \mu(u)}{H(u)} + \frac{\left\lfloor \frac{X(u)}{G(u)} \right\rfloor \cdot \phi(u) + \left\lfloor \frac{G(u) \cdot H(u)}{P(u)} \right\rfloor \cdot \delta(u) + \delta(u) \cdot \phi(u)}{H(u)} \right\rfloor, \tag{26}$$

where deg(δ($u$)), deg(ϕ($u$)) ≤ –1. We wish the second term in the above summation to have a degree less than 0. To achieve that we require

(A) $\deg\left(\left\lfloor \frac{X(u)}{G(u)} \right\rfloor\right) \le \deg(H(u))$;

(B) $\deg\left(\left\lfloor \frac{G(u) \cdot H(u)}{P(u)} \right\rfloor\right) \le \deg(H(u))$.

In the following, we assume these conditions to be satisfied. Thus (26) simplifies to

$$Q(u) = \left\lfloor \frac{\left\lfloor \frac{X(u)}{G(u)} \right\rfloor \cdot \mu(u)}{H(u)} + \Delta(u) \right\rfloor = \left\lfloor \frac{\left\lfloor \frac{X(u)}{G(u)} \right\rfloor \cdot \mu(u)}{H(u)} \right\rfloor, \tag{27}$$

as deg(Δ($u$)) ≤ –1. We also note that deg(μ($u$)) = deg($G$($u$)) + deg($H$($u$)) – deg($P$($u$)).

The above analysis leads to the conclusion summarized in the following theorem:

**Theorem 1.** Let $A$($u$), $B$($u$) and $P$($u$) be given such that 0 ≤ deg($A$($u$)), deg($B$($u$)) < deg($P$($u$)) = $N$. For the modular polynomial computation $X$($u$) mod $P$($u$), $X$($u$) = $A$($u$) · $B$($u$), if polynomials $G$($u$) and $H$($u$), α = deg($G$($u$)) and β = deg($H$($u$)), satisfy the conditions

$$\deg(X(u)) \le 2N - 2 \le \alpha + \beta;\ \alpha \le \deg(P(u)) = N, \tag{28}$$

then the quotient polynomial in (27) is same as the quotient polynomial $\lfloor X(u) / P(u) \rfloor$.

A generalization of the formulation of polynomials $G$($u$) and $H$($u$) away from polynomials of the type $u^a$ is crucial here. A simple choice of degrees that satisfy (28) is α = β = $N$. In fact, $G$($u$) and $H$($u$) can even be identical.

All the above analysis results in the new BA-P summarized in the following.

**A New Barrett Algorithm for Polynomials for Computing $A(u) \cdot B(u) \bmod P(u)$ (BA-P)**

Input: $A(u)$, $B(u)$, $P(u)$, $G(u)$, $H(u)$; $0 \le \deg(A(u))$, $\deg(B(u)) < N$, $N = \deg(P(u))$, $\alpha = \deg(G(u))$, $\beta = \deg(H(u))$.

Output: $Q(u)$ (quotient polynomial when $A(u) \cdot B(u)$ is divided by $P(u)$)

**Step 0.** Pre-compute $\mu(u) = \lfloor G(u) \cdot H(u) / P(u) \rfloor$, $\deg(\mu(u)) = \alpha + \beta - N$ (one-time)

**Step 1.** Compute $X(u) = A(u) \cdot B(u)$, $\deg(X(u)) \le 2 \cdot N - 2$ (ordinary multiplication)

**Step 2.** Compute $D(u) = \lfloor X(u) / G(u) \rfloor$, $\deg(D(u)) \le 2 \cdot N - \alpha - 2$ (quotient)

**Step 3.** Compute $E(u) = D(u) \cdot \mu(u)$, $\deg(E(u)) \le N + \beta - 2$ (ordinary multiplication)

**Step 4.** Compute $Q(u) = \lfloor E(u) / H(u) \rfloor$, $\deg(Q(u)) \le N - 2$ (quotient)

Once $Q(u)$ is computed, the algorithm may proceed further to compute the remainder $X(u) \bmod P(u)$ as

**Step 5.** $C(u) = X(u) - Q(u) \cdot P(u)$, $\deg(C(u)) \le N - 1$. (remainder, ordinary multiplication)

The conditions in (28) that are required for the polynomials $G(u)$ and $H(u)$ are general and leave door open to a wide range of possibilities. The corresponding computational steps can also be vastly different.

**Example 1.** Assume that the computation is defined in GF(2). Let $N = 6$, $P(u) = u^6 + u + 1$. We can choose $G(u) = u^6 + 1$, $H(u) = u^6 + 1$. Then $\mu(u) = u^6 + u + 1$. Let $X(u) = u^{10} + u^9 + u^8 + u^4 + u^2 + 1$. The various polynomials are $D(u) = u^4 + u^3 + u^2$, $E(u) = u^{10} + u^9 + u^8 + u^5 + u^2$. $Q(u) = u^4 + u^3 + u^2$, $Q(u) \cdot P(u) = u^{10} + u^9 + u^8 + u^5 + u^2$, $C(u) = u^5 + u^4 + 1$.

**Example 2.** Assume that the computation is defined in GF(2). An interesting case arises when $G(u) = H(u)$ and $\mu(u) = P(u)$. In that case, $G(u)^2 = P(u)^2 + R(u)$, $\deg(R(u)) < N$. If $G(u) = \sum_{i=0}^{N} G_i u^i$ and $P(u) = \sum_{i=0}^{N} P_i u^i$, then $G(u)^2 = \sum_{i=0}^{2N} G_i u^{2i}$ and $P(u)^2 = \sum_{i=0}^{2N} P_i u^{2i}$. A trivial solution to $G(u)^2 = P(u)^2 + R(u)$ is $G(u) = P(u)$, $R(u) =$ 0. Other interesting possibilities are $R(u) = \sum_{i=0}^{(N-1)/2} R_i u^{2i}$, where $R_i \in GF(2)$ and can take any values. This results in $G_i = P_i + R_i$, $i = 0, \ldots, (N-1)/2$; $G_i = P_i$, $i = (N+1)/2, \ldots, N$.

The analysis in Example 2, though applicable to polynomials defined in GF(2) only, can also be used to identify other desirable forms for $G(u)$ and $H(u)$ for a given $P(u)$ including those that require no pre-computations [19]. The analysis in [19] focuses only on $G(u) = u^N$ and $P_i = 0$ for $i = \lfloor N/2 \rfloor, \ldots, N-1$. In general, we have $G(u) \cdot H(u) = \mu(u) \cdot P(u) + R(u)$, $\deg(R(u)) < N$. So for a desirable form for $\mu(u)$ and arbitrary $R(u)$, $\deg(R(u)) < N$, numerous possibilities exist for the polynomials $G(u)$ and $H(u)$. Conversely, we can choose desirable forms for $G(u)$ and $H(u)$ that can lead to computational simplifications of steps 2 and 4. This will be the focus of our study in the next section.

We end this section by stating that the next section establishes RPS based BA-P in GF(2) though most results are general such that they are directly, or with suitable modifications, applicable to GF($p^a$), $p$ being a prime, and other number systems such as rational, real, and complex numbers.

## IV. RPS Based Barrett Algorithm for Polynomials

We now turn to a RPS based BA-P for the computation in Steps 1-5. All the polynomials need to be expressed in residue form in suitable RPS defined mod $M(u)$. Step 0 involves a one-time computation, so it is straightforward to map to RPS. Steps 1, 3, and 5 involve ordinary polynomial multiplication and hence they are also straightforward to map to RPS. Steps 2 and 4 require computation of quotient polynomial up on division by polynomials $G(u)$ and $H(u)$, respectively. As seen in Appendix A, in order to compute these two steps solely with polynomial residue arithmetic, both $G(u)$ and $H(u)$ must be a factor of $M(u)$ that defines the RPS used for computation. Also, since the BA-P is expected to be used recursively for carrying out MPE, we would like to use the same RPS in all Steps 1-5. In addition, for residues to correspond to the actual polynomials in Steps 1-5, $\deg(M(u))$ must exceed the maximum value of degree of polynomial at each step.

Given $P(u)$, $G(u)$, and $H(u)$ that satisfy Theorem 1, it is easy to calculate the smallest degree of RPS modulo polynomial $M(u)$ as it is at least one larger than the largest degree polynomial in steps 1-5. It is interesting to note that $\gcd(G(u), H(u)) = 1$ is not required. We are proposing to compute the quotients as required in steps 2 and 4 using the algorithm described in Appendix A. Since this algorithm requires modulo inverses, it works only when the various moduli polynomials are relatively co-prime. Also, the result of $\lfloor X(u) / G(u) \rfloor$ would only be known in terms of residues corresponding to the moduli that constitute $M(u) / G(u)$. Thus, for the quotient residues to be computed in step 2, we require $\gcd(G(u), M(u) / G(u)) = 1$. Similarly, it is also required that $\gcd(H(u), M(u) / H(u)) = 1$ for quotient residues in step 4.

Based on this discussion and the conditions in Theorem 1, we have the following additional conditions for the RPS based on $M(u)$ to be used in BA-P:

1. largest degree of polynomials in Steps 1-5 < $\deg(M(u)) = L$

2. $G(u) \mid M(u)$

4. $H(u) \mid M(u)$

5. $\gcd(G(u), M(u) / G(u)) = 1$

6. $\gcd(H(u), M(u) / H(u)) = 1$.

A large number of possibilities become apparent now. We can select $G(u)$ and $H(u)$ first that satisfy the above conditions. Then the polynomial $M(u)$ is constructed such that $\mathrm{lcm}(G(u), H(u)) \mid M(u)$. Finally, if needed further suitable residues are included in $M(u)$ to satisfy the first condition. Also it is possible to select $M(u)$ first and then select $G(u)$ and $H(u)$ in terms of factors of $M(u)$. For instance, if $\alpha = \beta = N$, then $L > 2 \cdot N - 2$.

The above description of RPS based BA-P via computation of quotient residues also brings out another aspect. The computations in steps 1, 3 and 5 are carried out mod $M(u)$. After Step 2 computation is done, the quotient $\lfloor X(u) / G(u) \rfloor$ is available in residues corresponding to $M(u) / G(u)$. Hence we need to carry out BEX-P to expand the quotient residues back to mod $M(u)$. Similarly, we need to carry out another BEX-P to expand the quotient residues of $\lfloor E(u) / H(u) \rfloor$ computed mod $M(u) / H(u)$ to mod $M(u)$. Such a BEX-P algorithm is described in Section II and also in Appendix B.

**A new RPS based Barrett Algorithm for Polynomials (BA-MPM)**

**Given:** $M(u)$, $G(u)$, $H(u)$. Let $M(u)$ have $n$ factors.

In step 2a, the first $a$ factors of $M(u)$ correspond to $G(u)$ while in step 4a, the first $b$ correspond to $H(u)$. This is assumed without any loss in generality.

**Input:** Residues of $A(u)$ and $B(u)$, that is, $(A_i(u), B_i(u) \equiv (A(u), B(u))$ (mod $M_i(u)$), $i = 1, \dots, n$.

**Pre-computational Step**

**Step 0.** Compute $\mu_i(u)$, $i = 1, \dots, n$, residues of $\mu(u) = \lfloor G(u) \cdot H(u) / P(u) \rfloor$.

**Computational Steps**

**Step 1. Modulo multiplication.** Compute $X_i(u) \equiv A_i(u) \cdot B_i(u)$, $i = 1, \dots, n$.

**Step 2a. Quotient computation.** Compute quotient residues $D_i(u)$, $i = a + 1, \dots, n$, from residues $X_i(u)$ , $i = 1, \dots, n$, and moduli $G_i(u)$, $i = 1, \dots, a$.

**Step 2b. BEX-P.** Use BEX-P on residues $D_i(u)$, $i = a + 1, \dots, n$, to get $a$ residues $D_i(u)$, $i = 1, \dots, a$.

**Step 3. Modulo multiplication.** Compute $E_i(u) \equiv D_i(u) \cdot \mu_i(u)$, $i = 1, \dots, n$.

**Step 4a. Quotient computation.** Compute quotient residues $Q_i(u)$, $i = b + 1, \dots, n$, from residues $E_i(u)$, $i = 1, \dots, n$, and moduli $H_i(u)$, $i = 1, \dots, b$.

**Step 4b. BEX-P.** Use BEX-P on residues $Q_i(u)$, $i = b + 1, \dots, n$, to get $b$ residues $Q_i(u)$, $i = 1, \dots, b$.

**Step 5. Remainder computation.** Compute $C_i(u) \equiv X_i(u) - Q_i(u) \cdot P_i(u)$, $i = 1, \dots, n$.

**END**

Provided that $L - \beta > N - 1$ (a rather trivial condition at this stage), Steps 4b and 5 may also be swapped to save $b$ MADD and MMULT. In such a case, we have the following steps:

**Step 5. Remainder computation.** Compute $C_i(u) \equiv X_i(u) - Q_i(u) \cdot P_i(u)$, $i = b + 1, \dots, n$.

**Step 6. BEX-P.** Use BEX-P on residues $C_i(u)$, $i = b + 1, \dots, n$, to get $b$ residues $C_i(u)$, $i = 1, \dots, b$.

**Example 3.** Assume that the computation is defined in GF(2). Let $N = 2^P - 1$, and $\alpha = \beta = N$. In this case, $\deg(M(u)) > 2 \cdot N - 2$. Let us choose $M(u) = (u^N - 1) \cdot [(u^{N+2} - 1) / (u - 1)]$ with $G(u) = H(u) = u^N - 1$. Also, $\gcd(u^a - 1, u^b - 1) = u^{\gcd(a,\, b)} - 1$. In our case, $N$ and $N + 2$ are two consecutive odd integers, hence $\gcd(N, N + 2) = 1$. Thus, $u^N - 1$ and $[(u^{N+2} - 1) / (u - 1)]$ are co-prime. The factorization of $u^N - 1$ in GF(2) using GF($2^P$) is well known with each factor having degree $P$ or less. Factorization of $u^{N+2} - 1$ can be obtained from the factorization of $u^A - 1$, $A = 2^{2 \cdot P} - 1$ as $2^{2 \cdot P} - 1 = (2^P - 1) \cdot (2^P + 1)$. Consequently, the factorization of $u^{N+2} - 1$ in GF(2) using GF($2^{2 \cdot P}$) can be obtained with each factor having degree $2 \cdot P$ or less. The entire computation $A(u) \cdot B(u) \bmod P(u)$ with $N = 2^P - 1$ can be performed using residue arithmetic where one-half the moduli have degree $P$ or less and the other half have degree $2 \cdot P$ or less. These techniques are easily extended to cover values of $N$ that are factors of $2^P - 1$. The various computations correspond to computing circular convolutions, a topic that has been studied extensively in technical literature.

**Example 4.** We pursue previous example with $P = 10$. Thus $N = \alpha = \beta = 1{,}023$. This is an example with a large value of $N$. In this case, $u^{1,023} - 1$ has factors of degree 10 or less. Similarly, $u^{1,025} - 1$ has factors of degree 20 or less. Let us choose $M(u) = (u^{1,023} - 1) \cdot [(u^{1,025} - 1) / (u - 1)]$ with $G(u) = H(u) = u^{1,023} - 1$. Here, we are assuming that $\deg(A(u)) = \deg(B(u)) = N - 1 = 1{,}022$. If the actual polynomials $A(u)$ and $B(u)$ as specified have degrees less than 1,022, then they may be padded with 0's and then treated as polynomials of degree 1,022. The entire computation $A(u) \cdot B(u) \bmod P(u)$ with $N = 1{,}023$ can be performed using residue arithmetic where one-half the moduli have degree 10 or less and the other half have degree 20 or less.

**Example 5.** Consider the case of the computation $A(u) \cdot B(u) \bmod P(u)$ defined in the field of real or complex numbers. We let $G(u) = H(u) = u^N - 1$ and $M(u) = u^{2 \cdot N} - 1 = (u^N - 1) \cdot (u^N + 1)$. Let $\omega$ be the $2 \cdot N$ root of unity. In this case, the entire computation can be carried out using DFT and IDFT via a suitably chosen FFT. For instance, $X(u)$ in step 1 can be computed using a size $2 \cdot N$ DFT. In step 2a, let $R(u)$ denote the remainder $X(u) \bmod G(u)$. Then $R(u)$ can be computed as a size $N$ IDFT of the even DFT coefficients of $X(u)$. We may write

$$D(u) = [X(u) - R(u)] / (u^N - 1).$$

Substituting odd powers of $\omega$, we get,

$$D(\omega^{2 \cdot k+1}) = -0.5 \cdot [X(\omega^{2 \cdot k+1}) - R(\omega^{2 \cdot k+1})],\ k = 0, \ldots, N - 1.$$

Computation of $R(\omega^{2 \cdot k+1})$ is same as a size $N$ DFT of sequence $R_i \cdot \omega^i$, $i = 0, \ldots, N - 1$. $D(\omega^{2 \cdot k+1})$ is same as a size $N$ DFT of sequence $D_i \cdot \omega^i$, $i = 0, \ldots, N - 1$. The computation of $D(u)$ in step 2a is performed by first taking size $N$ IDFT of $D(\omega^{2 \cdot k+1})$, $k = 0, \ldots, N - 1$, obtaining the sequence $D_i \cdot \omega^i$, $i = 0, \ldots, N - 1$, and then constructing $D_i$, $i = 0, \ldots, N - 1$, from it. The BEX-P in step 2b consists in taking size $N$ DFT of $D(u)$. Step 4 is similar to step 2. A total of 2 IDFT and 2 DFT, each of size $N$, are required to be computed in step 2. Same goes for step 4. Steps 0, 1, 3, and 5 are straightforward as will be clear to the reader from the above description.

We have not included a thorough computational complexity analysis here. The work is focused on residue polynomial arithmetic and algorithms. Based on our understanding of the various methodologies, we conjecture that the work will lead to implementations that are significantly faster in performing the computations underlying the new algorithms described here.

## V. Further Analysis

In this section, we describe an algorithm for MPE, called BA-MPE, that uses the new RPS based BA-MPM. For clarity, let the BA-MPM that computes $C(u) = A(u) \cdot B(u) \pmod{P(u)}$ be denoted by BA-MPM($A$, $B$).

**BA-MPE**

**Input:** Residue vector **A(*u*)** for $A(u)$, $P(u)$, and $E$, $E = \sum_{i=0}^{k} e_i \cdot 2^i$ .

**Output:** Residue vector **C(*u*)** for $C(u)$, $C(u) = A(u)^E \bmod P(u)$.

1. If $e_0 = 1$

   **C(*u*)** $\leftarrow$ **A(*u*)**

   else

   **C(*u*)** $\leftarrow$ **<u>1</u>**

2. For $j = 1$ to $k$ do

   **A(*u*)** $\leftarrow$ BA-MPM($A$, $A$)

      If $e_j = 1$ then

      **C(*u*)** $\leftarrow$ BA-MPM($C$, $A$)

   end If

   end For.

Here, **<u>1</u>** denotes the vector of all 1s. No further explanation is required for the BA-MPE.

## VI. Conclusions

In this work, new Barrett algorithms are described for computing modulo multiplication $A(u) \cdot B(u) \bmod P(u)$ and modulo exponentiation $A(u)^E \bmod P(u)$, $P(u)$ being an irreducible polynomial of degree $N$. A residue polynomial system based version of the new Barrett algorithm is also described that uses only residue arithmetic thereby avoiding large degree polynomial multiplication that may be computationally intensive. All the algorithms as described here are a first from our understanding of the research literature. The previously known Barrett algorithms use powers of $u$ to scale the various computations.

## APPENDIX A: Computing Quotient Residue Polynomials in RPS

### Problem Formulation.

Given residues $X_i(u)$ of a polynomial $X(u)$, $X_i(u) \equiv X(u)$ (mod $M_i(u)$), $i = 1, \ldots, n$,

$$M(u) = M_{\mathrm{I}}(u) \cdot M_{\mathrm{II}}(u),$$

$$M_{\mathrm{I}}(u) = \prod_{i=1}^{a} M_i(u),$$

$$M_{\mathrm{II}}(u) = \prod_{i=a+1}^{n} M_i(u),$$

$$\gcd(M_{\mathrm{I}}(u), M_{\mathrm{II}}(u)) = 1;$$

Compute residues of quotient $Q(u)$ when $X(u)$ is divided by $M_{\mathrm{I}}(u)$, $0 \le \deg(Q(u)) < \deg(M_{\mathrm{II}}(u))$.

We now revisit polynomial arithmetic introduced in Section II and use it to describe the algorithm for computing residues of $Q(u)$. Consider (1) when $X(u)$, $A_1(u)$ and $R_1(u)$ are known. We may compute $Q_1(u)$ as

$$Q_1(u) = (X(u) - R_1(u)) \cdot A_1(u)^{-1}. \qquad \text{(A1)}$$

Again, when $Q_1(u)$, $A_2(u)$ and $R_2(u)$ are known we may compute $Q_2(u)$ as

$$Q_2(u) = (Q_1(u) - R_2(u)) \cdot A_2(u)^{-1}. \qquad \text{(A2)}$$

This can be carried out recursively to finally compute $Q_a(u)$ as

$$Q_a(u) = (Q_{a-1}(u) - R_a(u) \cdot A_a(u)^{-1}. \qquad \text{(A3)}$$

The representation of $X(u)$ in (5) is still valid and is reproduced below for completeness.

$$X(u) = Q_a(u) \cdot [A_a(u) \cdots A_1(u)] + [R_a(u) \cdot (A_{a-1}(u) \cdots A_1(u)) + \ldots + R_2(u) \cdot A_1(u) + R_1(u)]. \qquad \text{(A4)}$$

Let us apply the above polynomial arithmetic in (A1)-(A4) to the RPS defined mod $M(u)$. We are given moduli $M_i(u)$ and corresponding residues $X_i(u)$, $i = 1, \ldots, n$. We now set

$$A_i(u) = M_i(u). \qquad \text{(A5)}$$

Thus,

$$R_1(u) = X_1(u). \tag{A6}$$

This leads to,

$$Q_1(u) = (X(u) - X_1(u)) \cdot M_1^{-1}(u) \tag{A7}$$

Since $X(u)$ is expressed in terms of its residues and $M_1^{-1}(u)$ exists only mod $M_i(u)$, $i = 2, \ldots, n$, we compute residues of $Q_1(u)$ in (A7) by taking mod $M_i(u)$, $i = 2, \ldots, n$, of both sides. This results in

$$Q_{1,i}(u) \equiv (X_i(u) - X_1(u)) \cdot M_1(u)^{-1} \ (\text{mod}\ M_i(u)),\ i = 2, \ldots, n. \tag{A8}$$

Since, $\deg(Q_1(u)) = \deg(X(u)) - \deg(M_1(u)) < \deg(M(u)) - \deg(M_1(u)) = \sum_{i=2}^{n} \deg(M_i(u))$, $Q_1(u)$ is uniquely expressed in terms of its residues $Q_{1,i}(u)$, $i = 2, \ldots, n$. After the first iteration in (A1), it is seen that

$$R_2(u) \equiv Q_2(u) \ (\text{mod}\ M_2(u)) = Q_{1,2}(u). \tag{A9}$$

Again, expressing (A9) in residue form, we compute residues of $Q_2(u)$ by taking (mod $M_i(u)$), $i = 3, \ldots, n$, of both sides. This results in

$$Q_{2,i}(u) \equiv (Q_{1,i}(u) - Q_{1,2}(u)) \cdot M_2(u)^{-1} \ (\text{mod}\ M_i(u)),\ i = 3, \ldots, n. \tag{A10}$$

Again, $\deg(Q_2(u)) < \sum_{i=3}^{n} \deg(M_i(u))$. Thus, $Q_2(u)$ is uniquely expressed in terms of its residues $Q_{2,i}(u)$, $i = 3, \ldots, n$.

This process is carried out for $a$ iterations computing residues of $Q_k(u)$ or $Q_{k,i}(u)$, $i = k + 1, \ldots, n$, at the $k$-th iteration, $k = 1, \ldots, a$. Up on conclusion of $a$ iterations, we have (A4) as

$$X(u) = Q_a(u) \cdot (M_a(u) \cdots M_1(u)) + [R_a(u) \cdot (M_{a-1} \cdots M_1(u)) + \ldots + R_2(u) \cdot M_1(u) + R_1(u)]. \tag{A11}$$

Thus, $Q_a(u)$ is the remainder obtained by dividing $X(u)$ by $M_A(u)$ expressed in terms of its residues $Q_{a,i}(u)$, $i = a + 1, \ldots, n$.

**A new algorithm to compute quotient polynomials**

**Input:**

RPS defined mod $M(u)$; $M(u) = M_{\mathrm{I}}(u) \cdot M_{\mathrm{II}}(u)$; residues of $X(u)$, $X_i(u) \equiv X(u) \pmod{M_i(u)}$, $i = 1, \ldots, n$.

**Output:**

Residues of $Q(u) \pmod{M_i(u)}$, $i = a + 1, \ldots, n$, $Q(u)$ being quotient polynomial when $X(u)$ is divided by $M_{\mathrm{A}}(u)$.

**Initialization:**

$Q_{0,i}(u) = X_i(u)$, $i = 1, \ldots, n$.

**Computational Steps:** For $k = 1, \ldots, a$, compute

$$Q_{k,i}(u) \equiv (Q_{k-1,i}(u) - Q_{k-1,k}(u)) \cdot M_k(u)^{-1} \pmod{M_i(u)},\ i = k + 1, \ldots, n.$$

**Output:**

$Q_i(u) = Q_{a,i}(u)$, $i = a + 1, \ldots, n$.

**END**

The first iteration ($i = 2, \ldots, n$) requires $n - 1$ MADD and MMULT. In general, $k$-th iteration ($i = k + 1, \ldots, n$) requires $n - k$ MADD and MMULT. The computations in each iteration can be performed in parallel in one computational step. There is a total of $a$ iterations requiring $(n - 0.5 \cdot a - 0.5) \cdot a$ MADD and MMULT. All the modular inverses are assumed to be pre-computed and stored.

## APPENDIX B: Computing Base Extension in RPS

**Problem Formulation.** BEX-P computational problem is:

Given

1. Two RPS defined mod $M_{\mathrm{I}}(u)$ and $M_{\mathrm{II}}(u)$.

2. Residues of $X(u)$ mod $M_I(u)$, that is, $X_i(u) \equiv X(u) \pmod{M_i(u)}$, $i = 1, \ldots, a$.

Compute

$n - a$ residues of $X(u)$ mod $M_{II}(u)$, that is, $X_i(u) \equiv X(u) \pmod{M_i(u)}$, $i = a + 1, \ldots, n$.

We wish not to involve CRT-P in BEX-P computation either via (15) or (17).

The representation of $X(u)$ in (5) is still valid and is reproduced below for completeness.

$$X(u) = Q_n(u) \cdot [A_n(u) \cdots A_1(u)] + [R_n(u) \cdot (A_{n-1}(u) \cdots A_1(u)) + \ldots + R_2(u) \cdot A_1(u) + R_1(u)]. \quad \text{(B1)}$$

Let us apply the above polynomial residue arithmetic to the RPS defined mod $M_A(u)$. We are given moduli $M_i(u)$ and corresponding residues $X_i(u)$, $i = 1, \ldots, a$. We now extend the moduli set to include moduli $M_i(u)$, $i = a + 1, \ldots, n$. Since the corresponding residues are unknown (and need to be computed), we extend the residue set to include residues $X_i(u)$, $i = a + 1, \ldots, n$, by choosing any arbitrary values (can be all 0s) for them. In the current formulation, we set

$$A_i(u) = M_i(u). \quad \text{(B2)}$$

Thus,

$$R_1(u) = X_1(u). \quad \text{(B3)}$$

This leads to,

$$Q_1(u) = (X(u) - X_1(u)) \cdot M_1(u)^{-1} \quad \text{(B4)}$$

Since $X(u)$ is expressed in terms of its residues and $M_1(u)^{-1}$ exists only $\pmod{M_i(u)}$, $i = 2, \ldots, n$, we compute residues of $Q_1(u)$ in (B4) by taking mod $M_i(u)$, $i = 2, \ldots, n$, of both sides. This results in

$$Q_{1,i}(u) \equiv (X_i(u) - X_1(u)) \cdot M_1(u)^{-1} \pmod{M_i(u)}, \; i = 2, \ldots, n. \quad \text{(B5)}$$

Since, $\deg(Q_1(u)) = \deg(X(u)) - \deg(M_1(u)) < \deg(M(u)) - \deg(M_1(u)) = \sum_{i=2}^{n} \deg(M_i(u))$, it is uniquely expressed in terms of its residues $Q_{1,i}(u)$, $i = 2, \ldots, n$. After the first iteration in (B5), it is seen that

$$R_2(u) \equiv Q_2(u) \text{ (mod } M_2(u)) = Q_{1,2}(u). \tag{B6}$$

Again, expressing (B6) in residue form, we compute residues of $Q_2(u)$ by taking (mod $M_i(u)$), $i = 3, \ldots, n$, of both sides. This results in

$$Q_{2,i}(u) \equiv (Q_{1,i}(u) - Q_{1,2}(u)) \cdot M_2(u)^{-1} \text{ (mod } M_i(u)),\ i = 3, \ldots, n. \tag{B7}$$

Again, $\deg(Q_2(u)) < \sum_{i=3}^{n} \deg(M_i(u))$ and $Q_2(u)$ is uniquely expressed in terms of its residues $Q_{2,i}(u)$, $i = 3, \ldots, n$.

This process is carried out for $a$ iterations computing residues of $Q_k(u)$ or $Q_{k,i}(u)$, $i = k + 1, \ldots, n$ at the $k$-th iteration, $k = 1, \ldots, a$. Up on conclusion of $a$ iterations, we have

$$X(u) = Q_a(u) \cdot (M_a(u) \cdots M_1(u)) + [R_a(u) \cdot (M_{a-1}(u) \cdots M_1(u)) + \ldots + R_2(u) \cdot M_1(u) + R_1(u)]. \tag{B8}$$

The second term in (B8) is the MRS digits of $X(u)$ corresponding to residues $X_i(u)$, $i = 1, \ldots, a$. Also, residues of $Q_a(u)$ are known for $M_i(u)$, $i = a + 1, \ldots, n$. Thus we have

$$\text{MRS digits of } X(u) = X(u) - Q_a(u) \cdot (M_a(u) \cdots M_1(u)). \tag{B9}$$

Expressing (B9) in residue form, we get residues in RPS mod $M_i(u)$ as

$$X_i(u) \equiv (X(u) \text{ (mod } M_i(u)) \equiv X_i(u) - Q_{a,i}(u) \cdot M_A(u) \text{ (mod } M_i(u)),\ i = a + 1, \ldots, n. \tag{B10}$$

Thus, BEX-P computes residues as per (B10). The Steps of the BEX-P algorithm are as follows:

**The new BEX-P algorithm**

**Input:**

RPS defined mod $M_I(u)$ and $M_{II}(u)$; residues of $X(u)$ mod $M_I(u)$, $X_i(u) \equiv X(u)$ (mod $M_i(u)$), $i = 1, \ldots, a$.

**Output:**

Residues of $X(u)$ in $M_{II}(u)$, $X(u)$ (mod $M_i(u)$), $i = a + 1, \ldots, n$.

**Initialization:**

Assign arbitrary values (can be all 0s) to $X_i(u)$, $i = a + 1, \dots, n$; and initialize $Q_0(u)$ as

$$Q_{0,i}(u) = X_i(u),\ i = 1, \dots, n. \tag{B11}$$

**Computational Steps:** For $k = 1, \dots, a$, compute

$$Q_{k,i}(u) \equiv (Q_{k-1,i}(u) - Q_{k-1,k}(u)) \cdot M_k(u)^{-1} \pmod{M_i(u)},\ i = k + 1, \dots, n. \tag{B12}$$

**Final Computational Step:** Compute the BEX-P residues as

$$X_i(u) \equiv X(u) \pmod{M_i(u)} \equiv X_i(u) - Q_{n,i}(u) \cdot M_A(u) \pmod{M_i(u)},\ i = a + 1, \dots, n. \tag{B13}$$

**END**

In general, the $k$-th iteration ($i = k + 1, \dots, n$) requires $n - k$ MADD and MMULT. The computations in each iteration can be performed in parallel in one computational step. There is a total of $a$ iterations requiring $(n - 0.5 \cdot a - 1) \cdot a + (n - a)$ MADD and MMULT. All the modular inverses are assumed to be pre-computed and stored. It will be clear to the readers that the final residues computed by the new BEX-P algorithm are unique and initial arbitrary values assigned to $X_i(u)$, $i = a + 1, \dots, n$, have no bearing on them. This is stated as a theorem in the following.

**Theorem 2.** The residue values computed via the BEX-P algorithm results in correct residues for any arbitrarily assigned initial values to $X_i(u)$, $i = a + 1, \dots, n$.

The BEX-P algorithm in this appendix needs to be compared to the one in Section II from the standpoint of implementation in a parallel and distributed manner.